\documentclass[letterpaper]{amsart} 
 
\usepackage{amssymb,latexsym} 
\usepackage{amsmath}
\usepackage[all]{xy}

\usepackage{psfrag}
\usepackage{epsfig}
 
\numberwithin{equation}{section}

\newtheorem{theorem}{Theorem}[section] 
\newtheorem{proposition}[theorem]{Proposition} 
 
\newtheorem{corollary}[theorem]{Corollary} 
\newtheorem{lemma}[theorem]{Lemma} 
 
\theoremstyle{definition} 
 
\newtheorem{remark}[theorem]{Remark}

\def\proof{\smallskip\noindent {\bf Proof.\ }}

\def\ZZ{\mathbb{Z}}

\def\FF{\mathbb{F}}

\begin{document}

\title{Congruence invariants of matrix mutation}

\author{Ahmet I. Seven}

\author{\.Ibrahim \"Unal}

\address{Middle East Technical University, Department of Mathematics, 06800, Ankara, Turkey}
\email{aseven@metu.edu.tr}

\email{iunal@metu.edu.tr}


\date{October 18, 2024}



\begin{abstract}

Motivated by the recent work of R. Casals on binary invariants for matrix mutation, we study the matrix congruence relation on quasi-Cartan matrices. 
In particular, we obtain a classification and determine normal forms modulo 4. As an application, we obtain new mutation invariants, which include the one obtained by R. Casals.



\end{abstract}



\maketitle

\section{Introduction}
\label{sec:intro}

Mutation is a well-known operation which defines an equivalence relation, called mutation equivalence, on skew-symmetric integer matrices. It is important for the applications to determine algebraic properties of matrices which are invariant under the mutation equivalence. For example, matrix congruence class of a skew-symmetric matrix is a mutation-invariant property, so any invariant of congruence is also a mutation invariant. 
Another basic idea to find mutation invariants is the notion of a quasi-Cartan matrix, which gives a class of \emph{symmetric} matrices naturally associated to a skew-symmetric matrix \cite{BGZ}. 
Therefore, it is a natural problem to classify symmetric quasi-Cartan matrices under the matrix-congruence relation and establish their relation to the mutation equivalence.
In this paper, we obtain a complete solution to this problem over $ \ZZ /  4\ZZ$ in a more general setup: we classify symmetric integer matrices whose diagonal elements are even numbers and determine normal forms using linear algebraic data. 
Furthermore, we determine the congruence classes that correspond to a mutation equivalence class. In particular, we obtain new mutation invariants, which include the one obtained recently in \cite{C}. We also answer questions that were left open in \cite{C}.

To state our results, we need some terminology. In what follows, by a matrix, we always mean a square integer matrix. For a matrix $ A $, by ``$ A $ modulo n'', we mean the matrix $ A \pmod{n}$, i.e. the matrix whose entries are the modulo $ n $ reductions of the entries of $ A $. We denote by $ diag(m_1A_1, m_2A_2,...,m_kA_k) $ the block diagonal matrix whose diagonal blocks are $ m_1$ copies of $ A_1 $, $ m_2 $ copies of $ A_2 $,..., $ m_k $ copies of $ A_k$, where $ A_i $'s are square matrices. 
We say that two matrices $ A $ and $ A' $ are congruent over $ \ZZ $ if there is an (integer) matrix $ P $ with $\det(P)=\mp 1$ such that $P^TAP=A'$.  Our first result is the following classification theorem for the congruence relation modulo 4:


\begin{theorem}\label{th:mod4}
Suppose that $ A $ is a symmetric $n\times n$ matrix whose diagonal entries are all even numbers. Then, there is a matrix $ P $ with $\det(P)=\mp 1$ such that $P^TAP \pmod{4}$ is a block diagonal matrix of the form 
\begin{equation}\label{eq:normal}
\mathrm{diag}(r
\left( \begin{array}{cc}
    {2} & {1} \\
    {1} & {2}   \\
       \end{array} \right),\,
s\left( \begin{array}{cc}
    {0} & {1} \\
    {1} & {0}   \\
       \end{array} \right), \,
t\left( \begin{array}{cc}
    {0} & {2} \\
    {2} & {0}   \\
       \end{array} \right), \,
p\left( \begin{array}{c}
    {2}   \\
       \end{array} \right),\,
m\left( \begin{array}{c}
    {0}   \\
       \end{array} \right))
\end{equation}
for some non-negative integers $p,r,s,t,m$. (We refer to these matrices as \emph{normal forms} of $ A $).
\end{theorem}
\noindent
An example for the theorem could be found in Section \ref{sec:mut-proof}.

Our next result, Theorem~\ref{th:mod4-unique}, will determine the uniqueness properties of the parameters in the normal forms \eqref{eq:normal}. These properties depend on some linear algebraic data, so we first recall some basic facts from the theory of alternating and quadratic forms over $\FF_2=\ZZ /2\ZZ$: 


We denote $ V=\ZZ^n $ and $\overline{V} =\ZZ^n /2\ZZ^n$ for convenience. Following \cite[Sec. 2.12]{J}, a quadratic form $q$ is an $\FF_2$-valued function on $\overline{V}$
having the following property:
$$q(u+v)=q(u)+q(v)+\Omega(u,v), (\mathrm{for\:\: all} \: u,v \in \overline{V}),$$ 
where $\Omega:\overline{V}\times \overline{V} \to \FF_2$ is an alternating bilinear form.
It is clear that the quadratic form $q$ completely determines the 
associated bilinear form $\Omega$; conversely, the form $q$ is determined by its values on a basis and $\Omega$. Recall (see e.g. \cite{D}) that 
there exists a $symplectic$ basis $\{e_1,f_1,...,e_r,f_r,h_1,...,h_p\}$ in $\overline{V}$ such that
$\Omega(e_i,f_j)=\delta_{i,j}$, 
and the rest of the values of $\Omega$ are $0$;
here $\delta_{i,j}$ is the Kronecker symbol. Let us write 
$\overline{V} _0=\{x \in \overline{V}: \Omega(x,v)=0, \: \mathrm{for \: any} \: v\in  \overline{V}\}$, which is 
the kernel of the form $\Omega$ .
If $q({\overline{V}} _0)=\{0\}$, then the Arf invariant of $q$ is defined as 
$$Arf(q)=\sum q(e_i)q(f_i).$$
It is well known from the theory of quadratic forms that $Arf(q)$ is independent of the
choice of the symplectic basis \cite{D}.

Two quadratic forms $q$ and $q'$ on $\overline{V} $ are 
isomorphic if there is a
linear isomorphism $T:\overline{V} \to \overline{V} $ such that $q(T(x))=q'(x)$ for any $x\in \overline{V}  $.
Isomorphism classes of quadratic forms $\{q\}$ on $\overline{V}$  
are determined by their Arf invariants and 
their restrictions $\{q \vert _ {\overline{V}_0}\}$ \cite[Sec. 2.12]{J}.
More precisely, for fixed dimensions of $\overline{V}$ and $\overline{V}_0$, there exist at most 
three isomorphism classes of quadratic forms $\{q\}$ and each isomorphism class is determined 
by precisely one of the following:

\begin{itemize} 
\item[{\rm(i)}] 
 $q(\overline{V}_0)=0, Arf(q)=1$
\item[{\rm(ii)}] 
$q(\overline{V}_0)=0, Arf(q)=0$
\item[{\rm(iii)}] 
$q(\overline{V}_0)=\FF_2$
\end{itemize}

Let us now suppose that $ A $ is a matrix as in Theorem \ref{th:mod4}, i.e. a symmetric $n\times n$ matrix whose diagonal entries are all even numbers.
 We denote by $ (\, ,\, ) $ the associated symmetric bilinear form on $ \ZZ^{n}$, so $ (v,w)=v^{T}Aw $ for any two vectors $v,w $ in $ \ZZ^{n}$, where 
$ v^{T}$ is the transpose of $ v $. One may note that, since the diagonal entries of $ A $ are even,  for any $v$ in $ \ZZ^{n}$, the product $ (v,v) $ is an even number.
Then we define the quadratic form $q=q_A$ associated with $ A $  as the function defined on $ \overline{V} $  as $ q(v)=(v,v)/2 \pmod{2} $, for all $ v $ in $ \overline{V} $ (where, by an abuse of notation, we write $ v $ to denote its coset $ v+2\ZZ $). It follows from a direct check that $ q $ is a quadratic form whose associated bilinear form $ \Omega $ is defined as $\Omega(v,w)=(v,w) \pmod{2} $.  



Associated with $ A $, we consider the following $ \ZZ $-submodules of $ V=\ZZ^n $:

\begin{itemize} 
\item
$V_{0}=\{v \in V: (v,w) \equiv 0\pmod{2} \mathrm{ \, for \, all \, } w \mathrm{ \, in \,} V\} $,
\item 
$V_{00}=\{v \in V_{0} : (v,v) \equiv 0 \pmod{4} \} $,
\item
$V_{000}=\{v \in V_{0} : (v,w) \equiv 0 \pmod{4} \mathrm{ \, for \, all \, } w \mathrm{ \, in \,} V \} $.
\end{itemize}




\noindent
Let us note that $ V_{00} $ is a submodule of $ V_0 $ and $ V_{000} $ is a submodule of $ V_{00} $.
Then we have the following quotient modules:

\begin{equation}\label{eq:V000}
\overline{V}_{0}= V_{0}/2V, \quad  \overline{V}_{00}= V_{00}/2V, \quad   \overline{V}_{000}= (V_{000}+2V)/2V
\end{equation}

\noindent
It follows from a direct check that these are vector spaces over $ \FF_2 =\ZZ /2\ZZ$ such that $ \overline{V}_{00} $ is a subspace of $ \overline{V}_0 $ and $ \overline{V}_{000} $ is a subspace of $ \overline{V}_{00} $.
The subspace $ \overline{V}_{0} $ is the kernel of the bilinear form $ \Omega $ on $ \overline{V} $ and $\overline{V}_{00}=\overline{V}_0 \cap q^{-1}(0)=\{v \in \overline{V}_0: q(v)=0\} $.
In particular, since $ \Omega $ is determined by $ q $ up to isomorphism, the dimensions $ \dim(\overline{V}_0) $ and $ \dim(\overline{V}_{00}) $ are determined by the isomorphism class of $ q $. 


Our next result shows that the congruence class of $ A $ modulo $ 4 $ is determined by the dimension of $ \overline{V}_{000}  $ and the isomorphism class of the corresponding quadratic form $ q=q_A$. More precisely:






\begin{theorem}\label{th:mod4-unique}
Let $ A $ be a a symmetric $ n \times n $ matrix whose diagonal entries are all even numbers and $ q $ be the associated quadratic form.
Then, for fixed dimension of $ \overline{V}_{000}  $, the parameters  $ p,r,s,t,m$ in the normal form \eqref{eq:normal} are uniquely determined by the isomorphism class of $ q $ up to the following: 


\begin{itemize} 
\item
$ m=\dim(\overline{V}_{000})    $,
\item
$p \equiv \dim(\overline{V}_{0}/\overline{V}_{000}) \pmod{ 2} $, 

\item
$ r\equiv Arf(q) \pmod{ 2}$ (if $q(\overline{V}_0)=0$).


\end{itemize}

In particular, for the minimum non-negative values of $ p$ and $r $, one gets unique normal form. More explicitly, if $q(\overline{V}_0)=0$ and $r \in \{0,1\} $, then  $ r $ (and others) are uniquely determined; if $q(\overline{V}_0)=\FF_2$ and $ p \in \{1,2\} $, $ r=0 $, then $ p $ (and others) are uniquely determined.

\end{theorem}

\noindent
Note that, in the theorem, the isomorphism class of the quadratic form $ q $ determines $q(\overline{V}_0)$ and $ \dim (\overline{V}_{0})=2t+p+m=n-2(r+s) $. It also determines whether the parameter $ p $ is equal to zero or not: $ p=0 $ if and only if $q(\overline{V}_0)=0$ (Lemma \ref{lem:q-V0}(vi)). 

As a corollary of this result, we obtain the following characterization of the determinants of symmetric matrices:

\begin{corollary}\label{cor:det}

Let $ A $ be a a symmetric $n\times n$  matrix whose diagonal entries are all even numbers. Let $ \overline{A} $ denote $ A \pmod{2} $. Then we have the following:

\begin{itemize} 
\item
$ \det(A)=\pm 1 \pmod{4}$ if and only if rank of $ \overline{A} $ is $ n $ (i.e. $ \dim(\overline{V}_{0})=0 $),
\item
$  \det(A)=2 \pmod{4}$ if and only if $ n $ is odd such that the rank of $ \overline{A} $ is $ n-1 $ and $q(\overline{V}_0)=\FF_2$ (i.e. $\dim(\overline{V}_{0})=1, \dim(\overline{V}_{00})=0   $),
\item
$  \det(A)=0 \pmod{4}$ if and only if $ \dim(\overline{V}_{00}) \geq 1$.
\end{itemize} 
In particular, if the rank of $ \overline{A} $ is less than or equal to $ n-2 $, then $  \det(A) \equiv 0 \pmod{4}$.

\noindent
Furthermore, for $ n $ even, we have the following:
\begin{itemize} 
\item
if $ n\equiv 0 \pmod{4}$, then $ \det(A)= 1 \pmod{4}$ or $  \det(A)=0 \pmod{4}$
\item
if $ n\equiv 2 \pmod{4}$, then $ \det(A)= -1 \pmod{4}$ or $  \det(A)=0 \pmod{4}$
\end{itemize} 
according to whether the rank of $ \overline{A} $ is $ n $ or less than $ n $ respectively.

%



\end{corollary}

\begin{remark}\label{rem:det}
According to the corollary, for $ n $ even, $ \det(A) \pmod{4} $ is a binary function, i.e. it takes only two values. The same property holds if $ n $ is odd, because, then the rank of $ \overline{A} $ is less than or equal to $ n-1 $ , so $  \det(A) \equiv 0 \pmod{4}$ or $  \det(A) \equiv 2 \pmod{4}$ respectively. 

Let us also note that $ \det(A) \pmod{4}$ is determined by $ \dim(\overline{V}_{0}) $ and $ \dim(\overline{V}_{00}) $ as in the corollary. However, if $ n $ is even, then it is determined only by $ \dim(\overline{V}_{0}) $.

This implies the following: if $ A $ and $ A' $ are as in the corollary such that $ \overline{A} $ and $ \overline{A'} $ have the same rank and the associated quadratic forms $ q_A $ and $ q_{A'} $ are isomorphic, then $  \det(A) \equiv \det(A') \pmod{4}$. (Note that $\overline{V}_{00}=\overline{V}_0 \cap q^{-1}(0)=\{v \in \overline{V}_0: q(v)=0\} $).

\end{remark}


An important class of matrices that satisfy the hypothesis of Theorem \ref{th:mod4} are symmetric quasi-Cartan matrices \cite{BGZ}. More precisely, a symmetric $ n \times n $ matrix $A$ is a quasi-Cartan matrix if and only if $A_{i,i}=2$ for all $i=1,2,..,n $. These include the well-known (symmetric) generalized Cartan matrices of Lie theory \cite{K}. In cluster algebra theory, those matrices are considered as companions of skew-symmetric matrices: a \emph{quasi-Cartan companion} of a skew-symmetric matrix $B$ is a quasi-Cartan matrix $A$ such that $|A_{i,j}|= |B_{i,j}|$ for all $i \ne j$. This notion was introduced in \cite{BGZ} to study properties of the matrix mutation operation. 

To be more precise we briefly recall the definition of matrix mutation. This is an operation defined on skew-symmetrizable matrices, which is a class of matrices that include skew-symmetric matrices. 
\emph{In this paper, we only work with skew-symmetric matrices}, so we don't recall the definition of a skew-symmetrizable matrix. 
For each matrix index $ k $, the mutation in direction $ k $ transforms a skew-symmetric matrix $ B $ into 
$B'=\mu_k(B)$ whose entries are given as follows, where we use the notation $[b]_+ = \max(b,0)$ (see, e.g., \cite[Equation~(1.1)]{BGZ}):

\begin{equation}
\label{eq:matrix-mutation}
B'_{ij} =
\begin{cases}
-B_{ij} & \text{if $i=k$ or $j=k$;} \\[.05in]
B_{ij} + [B_{ik}]_+ [B_{kj}]_+ - [-B_{ik}]_+ [-B_{kj}]_+
 & \text{otherwise.}
\end{cases}
\end{equation}
This operation is involutive, so it defines a \emph{mutation-equivalence} relation on skew-symmetric matrices: two matrices are called mutation equivalent if they can be obtained from each other by a sequence of mutations. As discussed in \cite[Section 1]{C}, it is a natural problem to find invariant properties of the matrix mutation. The most basic mutation invariant is the corresponding congruence class: if $ B $ and $ B' $ are mutation equivalent, then there is a matrix $ P $ with $\det(P)=\mp 1$ such that $P^TBP= B' $  (see, e.g., \cite[Section 1]{C}). 

Another basic mutation invariant, which is not discussed in \cite{C}, is the corresponding quadratic form \cite[Proposition 5.2]{S2}. To be more precise, let us 
suppose that $ B $ is a skew-symmetric matrix and let $ \overline{B} $ denote the matrix $ B \pmod{2} $. 
Let $ \Omega $ be the alternating bilinear form on $ \overline{V}= \ZZ^n/2\ZZ^n $  whose Gram matrix with respect to a basis $ \overline{\mathcal{B}}  $ of  $ \overline{V} $ is $ \overline{B} $. 
Then there is a unique quadratic form $ Q=Q_B $ defined as follows: $Q(u+v)=Q(u)+Q(v)+\Omega(u,v)$, for all $ u,v \in  \overline{V}$ and $Q(b)=1$ for all $b \in \overline{\mathcal{B}}$. 
Furthermore, the isomorphism class of $ Q $ is invariant under the mutation operation:

\begin{proposition}\label{prop:skew-quad}\cite[Proposition 5.2]{S2}, \cite[Proposition 8.4]{BHII}
Suppose that $ B  $ and $ B' $ are skew-symmetric matrices which are mutation equivalent. Then the corresponding quadratic forms $ Q_B $ and $ Q_{B'} $ are isomorphic.

\end{proposition}
On the other hand, for any quasi-Cartan companion $ A $ of $ B $, we have the associated quadratic form $ q=q_A $. Let us note that, if we take the same initial basis $ \mathcal{B} $ to define $Q= Q_B $ as the one used for $q= q_A $, then, for any $ b $ in $ \mathcal{B} $, we have the following:
$q(b)=1$ (because $ (b,b)=2 $) and $ Q(b)=1 $ (by definition), so $q(b)=Q(b)$ (where, by an abuse of notation, we write $ v $ to denote its coset $ v+2\ZZ $). Also the matrices $ B \pmod 2$ and $ A \pmod 2$ coincide. Thus, in this setup, the quadratic forms $ Q_B $ and $ q_A $ also coincide. 
Then we have the following corollary of the above proposition:

\begin{corollary}\label{cor:skew-quad} 
Suppose that $ B  $ and $ B' $ are skew-symmetric matrices which are mutation equivalent. Let $ A $ and $ A' $ be quasi-Cartan companions of $ B $ and $ B' $ respectively. Then the corresponding quadratic forms $ q_A $ and $ q_{A'} $ are isomorphic.

\end{corollary}

Our next result classifies quasi-Cartan companions of mutation equivalent matrices modulo $ 4 $ and show that, for any quasi-Cartan companion $ A $, the number $ \det{A}\pmod{4}  $ is a general mutation invariant:  

\begin{theorem}\label{th:mut}
Suppose that $ B $ and $ B' $ are skew-symmetric matrices which are mutation equivalent. Let $ A $ and $ A' $ be quasi-Cartan companions of $ B $ and $ B' $ respectively. Let $ m $ and $ m' $ be the corresponding parameters in the normal forms of $ A $ and $ A' $ as in Theorem \ref{th:mod4}.
Then there is a matrix $ P $ with $\det(P)=\mp 1$ such that $P^TAP\equiv A' \pmod{4}$  if and only if $ m = m' $.
Furthermore, $ \det{A}\equiv \det{A'} \pmod{4}  $.
\end{theorem}

\noindent
We give an example of companions with different parameters $ m $ and $ m' $ in Section \ref{sec:mut-proof}.

We now discuss applications of our results to a particular invariant, called $ \delta $ invariant, introduced by R. Casals in \cite[Definition 1.3.]{C}, where for any skew-symmetric matrix $ B $, a particular quasi-Cartan companion $\mathfrak{S}(B)$ is defined as

\vspace{0.3cm}

$\mathfrak{S}(B):=V(B)+V(B)^t$, where 
$V(B)_{i,j} := \begin{cases}
  b_{ij}  & \text{ if }i<j \\
  1  & \text{ if }i=j \\
  0  & \text{ if }i>j \\
\end{cases}
$

\vspace{0.3cm}

\noindent
and $\delta(B)$ is defined as the determinant of $\mathfrak{S}(B) $ modulo 4. 




The main result of \cite{C} is that $ \delta $ is a mutation invariant. Our Theorem \ref{th:mut} shows that mutation invariance is not special to this particularly defined $ \delta $: the determinant modulo $ 4 $ of \emph{any} quasi-Cartan companion $ A $ is invariant under the mutation operation. More importantly, our results provide far reaching generalizations of these properties: Theorems \ref{th:mod4} and \ref{th:mod4-unique} classify, among others, all quasi-Cartan matrices (not only their determinants), including the particularly defined $\mathfrak{S}(B) $, under the matrix congruence relation modulo $ 4 $ and Theorem \ref{th:mut} determines their mutation invariance (as quasi-Cartan matrices).  This more algebraic approach allows us, in particular, to answer all questions that were left open in \cite{C}. More precisely, Corollary \ref{cor:det} gives a complete description of the $ \delta $ invariant, completing the proof that $ \delta $ is binary valued for any size $ n $. Furthermore, Corollary \ref{cor:det} determines exactly when $ \delta  $ is not a congruence invariant. We collect these properties as a corollary:


\begin{corollary}\label{cor:delta}


For any quasi-Cartan companion $ A $ of $ B $, we have $  \delta(B) = \det(A) \pmod{4}$, so $ \delta(B) $ is determined as in Corollary \ref{cor:det}. In particular, $ \delta(B) $ is binary.

Suppose that $ B $ and $ B' $ are  skew-symmetric $ n \times n $ matrices
and let $ \overline{B} $ (resp. $ \overline{B'} $) be their modulo $ 2 $ reductions. Suppose also that 
$ \overline{B} $ and $ \overline{B'} $ have the same rank, say $ \ell $ (e.g. when $ B $ and $ B' $ are congruent over $ \ZZ $). Then $ \delta(B)\neq \delta(B')$ if and only if 
$ n $ is odd with $ \ell=n-1 $ 
and $ Q_B(\overline{V}_{0})\ne Q_{B'}(\overline{V}'_{0}) $ (where $ \overline{V}_{0} $ and $ \overline{V}'_{0} $ are the kernels of $ \overline{B} $ and $ \overline{B'} $ respectively).


\end{corollary}
\noindent
To illustrate these properties, we give an example in Section \ref{sec:delta-ex}. 

\section{Proofs}
\label{sec:proof}

For the proofs, we suppose that $ A $ is as in the statement of the theorem. Let us recall from Section \ref{sec:intro} that $ (\, ,\, ) $ denotes the symmetric bilinear form defined by $ A $ on $ \ZZ^{n}$, so $ (v,w)=v^{T}Aw $ for vectors $v,w $ in $ \ZZ^{n}$, where $ v^{T}$ is the transpose of $ v $. We say that $ A $ is the Gram matrix of a $ \ZZ $-basis $ \mathcal{B}=\{v_1,v_2,...,v_n\} $ if $ A_{i,j}=(v_i,v_j) $ for all $ i,j=1,..,n $.

\subsection{Proof of Theorem~\ref{th:mod4}}

To prove the theorem, we first note that, since the diagonal entries are even numbers, the matrix $ \overline{A}=A\pmod{2}$ defines an alternating bilinear form on the $\FF_2$-vector space $V=\ZZ^n /2\ZZ^n$ by $ (v,w) \pmod{2}  $. 
Then, it is well known, (see e.g. \cite[Ch. XV, Theorem 8.1]{L}), that there is an invertible matrix, say $ \overline{S} $, over the field $\FF_2$ such that $ \overline{S}^T\overline{A}\overline{S} $ is a block diagonal matrix of the form 
$$\mathrm{diag}(k
\left( \begin{array}{cc}
    {0} & {1} \\
    {1} & {0}   \\
       \end{array} \right),\,
\ell \left( \begin{array}{c}
    {0}   \\
       \end{array} \right)
).
$$
Here the matrix $ \overline{S} $ is a matrix of $ 0 $'s and $ 1 $'s and invertible over $ \FF_2 $ (so $ \det(\overline{S})=1 $ in $ \FF_2 $). On the other hand, any invertible matrix over $ \FF_2 $ is  a product 
of elementary matrices \cite[Ch. XIII, Proposition 9.1]{L}, so the matrix $ \overline{S} $ is a product of elementary matrices, say $ \overline{S}=E_1 \cdot E_2 \cdots E_j $. Let us note that each of these elementary matrices is upper or lower triangular with $ 1 $'s on the diagonal, so its determinant considered over integers is also equal to $ 1 $, so each $ E_i $ is also invertible over $ \ZZ $. Then the product $ S=E_1\cdot E_2 \cdots E_j $ of these elementary matrices over integers is also invertible over $ \ZZ $ (here $\overline{S}= S \pmod{2} $).
Thus there is a $ \ZZ $-basis $ \mathcal{E}=\{e_1,f_1,...,e_{k},f_{k},g_1,...,g_{\ell}\} $ whose Gram matrix is $ S^TAS $ . 
(Note that $ 2k $ is equal to the rank of the matrix $A\pmod 2$). 
We will prove the theorem by obtaining a basis $ \mathcal{B} $ from $  \mathcal{E} $ such that the Gram matrix modulo $ 4 $ of $ \mathcal{B} $ is as in the statement of the theorem. 

To proceed, let us first note the following properties of the symmetric bilinear form $ ( \,,\, ) $ on $ \ZZ^{n} $:

$ (v,w)\equiv 0 \pmod{2}$ if and only if $ (v,w)\equiv 0$ or $ 2\pmod{4}$,

$ (v,w)\equiv 1 \pmod{2}$ if and only if $ (v,w)\equiv \pm 1 \pmod{4}$.

\noindent
Then, for any $v,w \in \mathcal{E}$, we have the following: 

$(v,w)\equiv \pm 1  \pmod 4$ if and only if $\{v,w\}=\{e_i,f_i\}$, for some $i$; if $\{v,w\}\ne \{e_i,f_i\}$, then $(v,w) \in \{0,2\} \pmod 4$.

\noindent
For convenience, by applying sign changes if necessary, we will assume that $ (e_i,f_i) \equiv 1  \pmod 4$ for any $ i $ (note that $ 2\equiv -2 \pmod{4} $).





In this setup, we first prove the following weaker version of Theorem~\ref{th:mod4}:
\begin{lemma}
\label{lem:mod4}
There is a  matrix $ P $ with $\det(P)=\mp 1$ such that $P^TAP \pmod{4}$ is a block diagonal matrix of the form 
$$\mathrm{diag}(r
\left( \begin{array}{cc}
    {2} & {1} \\
    {1} & {2}   \\
       \end{array} \right),\,
s\left( \begin{array}{cc}
    {0} & {1} \\
    {1} & {0}   \\
       \end{array} \right), \,
p\left( \begin{array}{c}
    {2}   \\
       \end{array} \right),\,
Z),
$$
where $Z$ is a (symmetric) matrix whose diagonal entries are all $0$'s and any non-zero entry is equal to $2$,  with $p,r,s$ being non-negative integers.
\end{lemma}
\proof We prove the lemma by obtaining  a $\ZZ$-basis from $\mathcal{E}$ whose Gram matrix modulo $4$ is as in its statement. We do this by replacing, for each $ i $, those basis vectors $v$ with $ (v,e_i)\equiv 2 $ or $ (v,f_i)\equiv 2 \pmod 4$. To proceed, suppose that there is $v$ in $\mathcal{E}$ with $v \notin \{e_i,f_i\} $ such that $ (v,e_i)\equiv 2 $ or $ (v,f_i)\equiv 2 \pmod 4$ . We may assume, without loss of generality, that $ (v,e_i)\equiv 2 \pmod 4 $. Let $ \mathcal{E'}$ be the basis obtained from $ \mathcal{E}$ by replacing $v$ with $v'=v+2f_i$. This is an invertible operation over integers (with the inverse replacing $ v' $ by $ v'-2f_i $), so $ \mathcal{E'}$ is a $ \ZZ $-basis. Furthermore,
$(v',e_i)=(v+2f_i, e_i)=(v,e_i)+2(f_i,e_i)=2+2\cdot 1\equiv 0 \pmod 4$ and, for any $w$ in $\mathcal{E'}$ with $w \ne e_i $, we have  $(v',w)=(e_i+2f_i,w)=(e_i,w)+2(f_i,w)\equiv (e_i,w) \pmod 4$ because, for $w\ne e_i$, $(f_i,w) \in \{0,2\}$. Also $(v',v')=(v+2f_i,v+2f_i)=(v,v)+4\cdot(v,f_i)+4\cdot(f_i,f_i)\equiv (v,v) \pmod 4$.
Thus the Gram matrix modulo $4$ of $ \mathcal{E'}$ is obtained from the Gram matrix of $ \mathcal{E}$ by making the entries corresponding to the basis vectors $e_i$ and $v$ equal to $0$ and keeping the other entries intact. 
Applying this base change operation to \emph{all} $v$ in $\mathcal{E}$ such that $ (v,e_i)\equiv 2 $ or $ (v,f_i)\equiv 2 \pmod{4} $ for some $ i $, we obtain a $ \ZZ $-basis $ \mathcal{B'}=\{e'_1,f'_1,...,e'_{k},f'_{k},g'_1,...,g'_{\ell}\} $ whose Gram matrix modulo $4$ is block diagonal of the form 
$$\mathrm{diag}(r
\left( \begin{array}{cc}
    {2} & {1} \\
    {1} & {2}   \\
       \end{array} \right),\,
s\left( \begin{array}{cc}
    {0} & {1} \\
    {1} & {0}   \\
       \end{array} \right), \,
\,
Z'),
$$
where each $ 2 \times 2 $ block is the Gram matrix that corresponds to a (hyperbolic) pair $ \{e_i,f_i\} $ and  $Z'$  is the submatrix (Gram matrix) that corresponds to the basis vectors $ g'_1,g'_2,...,,g'_{\ell} $.

To complete the proof of the lemma, we will further simplify $ Z' $ by ``isolating'' inductively those $ g'_i $ such that  $ (g'_i,g'_i) \equiv 2 \pmod 4$. More explicitly, let us assume that there exists $ g'_i $ such that  $ (g'_i,g'_i) \equiv 2 \pmod 4$. Without loss of generality, we may assume $ i=1 $. Then, for any $ j\ne 1  $, we have the following: if $  (g'_1,g'_j) \equiv 2 \pmod 4 $, then $  (g'_1,g'_j+g'_1) =(g'_1,g'_j)+(g'_1,g'_1)=2+2 \equiv 0 \pmod 4$. 
Thus, replacing any such $ g'_j $ by $ g''_j= g'_1+g'_j$ gives a $ \ZZ $-basis such that the entries corresponding to the basis vectors $g'_1$ and $g''_j$ in the Gram matrix are $0 $ modulo $4$.
Applying such base changes for all $g'_j  $ with $  (g'_1,g'_j) \equiv 2 $, $ 1\ne j  $, we obtain a basis $ \{e'_1,f'_1,...,e'_{k},f'_{k},g'_1,g''_2...,g''_{\ell}\} $ whose Gram matrix modulo $4$ is of the form
$$\mathrm{diag}(r
\left( \begin{array}{cc}
    {2} & {1} \\
    {1} & {2}   \\
       \end{array} \right),\,
s\left( \begin{array}{cc}
    {0} & {1} \\
    {1} & {0}   \\
       \end{array} \right), \,
\left( \begin{array}{c}
    {2}   \\
       \end{array} \right),\,
Z'').
$$
Now applying the same type of base changes inductively for all $g''_i  $, $ i\geq 2 $, with $  (g''_i,g''_i) \equiv 2 \pmod{4}$, we obtain a basis whose Gram matrix modulo $4$ is as in the statement of the lemma. This completes the proof of the lemma.

To complete the proof of Theorem~\ref{th:mod4}, we will show that the matriz $ Z $ in Lemma \ref{lem:mod4} is congruent to a block diagonal matrix:
\begin{lemma}
\label{lem:mod4-2}
Suppose that $Z$ is a symmetric matrix such that for any diagonal entry $ Z_{ii} $, we have $ Z_{ii} \equiv 0 \pmod{4}$ and, for any non-diagonal entry $ Z_{ij} $, $ i\ne j $, we have $ Z_{ij} \equiv 0 \pmod{4}$ or $ Z_{ij} \equiv 2 \pmod{4}$. 
Then there is a matrix $ P $ with $\det(P)=\mp 1$ such that $P^TAP\pmod{4}$ is a block diagonal matrix of the form 
$$\mathrm{diag}(t
\left( \begin{array}{cc}
    {0} & {2} \\
    {2} & {0}   \\
       \end{array} \right),\,
m\left( \begin{array}{c}
    {0}   \\
       \end{array} \right)
)
$$
for some non-negative integer $m$.
\end{lemma}

\proof We prove the lemma by an argument which is similar to the one at the beginning of this section. Let us write $ Z=2Z' $ and denote $ \overline{Z'}= Z' \pmod{2} $. Considering  $ \overline{Z'} $ as a matrix over $ \FF_2 $,
there is an invertible matrix $ \overline{P} $ over $ \FF_2 $ such that $ \overline{P}^T\overline{Z'}\overline{P}$ (over $ \FF_2 $)
is a block diagonal matrix of the form 
$$\mathrm{diag}(t
\left( \begin{array}{cc}
    {0} & {1} \\
    {1} & {0}   \\
       \end{array} \right),\,
m\left( \begin{array}{c}
    {0}   \\
       \end{array} \right)
)
$$
for some non-negative integer $m$ (in fact, $m$ is equal to the corank of $\overline{Z'}$). Here the matrix $ \overline{P} $ is a matrix of $ 0 $'s and $ 1 $'s and $ \det(\overline{P})=1$ in $ \FF_2 $. 
On the other hand, any invertible matrix over $ \FF_2 $ is a product 
of elementary matrices \cite[Proposition 9.1]{L}; in particular, $ \overline{P} $ is a product of elementary matrices, say $ \overline{P}=E_1\cdot E_2 \cdots E_j $. Each of these elementary matrices is upper or lower triangular with $ 1 $'s on the diagonal, so its determinant considered over integers is also equal to $ 1 $, so each $ E_i $ is also invertible over $ \ZZ $. Then the product $ P=E_1\cdot E_2 \cdots E_j  $ of these elementary matrices over integers is also invertible over $ \ZZ $; here $\overline{P}\equiv P \pmod{2}$. 
Let $ Z''= P^TZ'P$. Then $ Z'' \equiv \overline{P}^T\overline{Z'}\overline{P} \pmod{2}  $, so $Z'' \pmod{2} $ is block diagonal as above. Then all entries in the matrix $ Z''$ are even except those entries which are equal to $1$ in the block diagonal matriz $Z'' \pmod{2} $. This implies that $2Z''$ modulo $4$
is equal to the block diagonal matrix  
$$\mathrm{diag}(t
\left( \begin{array}{cc}
    {0} & {2} \\
    {2} & {0}   \\
       \end{array} \right),\,
m\left( \begin{array}{c}
    {0}   \\
       \end{array} \right)
)
$$
(any even number becomes congruent to $0$ modulo $4$ when multiplied by $2$). Thus $ P^TZP=P^T(2Z')P=2P^TZ'P=2Z''$ modulo $ 4 $ is as in the statement of the lemma. This completes the proof of the lemma and the proof of Theorem~\ref{th:mod4}

\subsection{Proof of Theorem~\ref{th:mod4-unique}}
We will prove the theorem after establishing some uniqueness properties of the associated parameters. We first show the following property of the parameter $ p $:
\begin{lemma}\label{lem:p}
Suppose $ A $ is a matrix such that $ A \pmod{4}$ is the matrix 
$$\mathrm{diag}(r
\left( \begin{array}{cc}
    {2} & {1} \\
    {1} & {2}   \\
       \end{array} \right),\,
s\left( \begin{array}{cc}
    {0} & {1} \\
    {1} & {0}   \\
       \end{array} \right), \,
t\left( \begin{array}{cc}
    {0} & {2} \\
    {2} & {0}   \\
       \end{array} \right), \,
p\left( \begin{array}{c}
    {2}   \\
       \end{array} \right),\,
m\left( \begin{array}{c}
    {0}   \\
       \end{array} \right))
$$
for some non-negative integers $r,s,t,m$ and $ p \geq 3 $.
Then, there is a matrix $ P $ with $\det(P)=\mp 1$ such that $P^TAP \pmod{4}$ is
$$\mathrm{diag}(r
\left( \begin{array}{cc}
    {2} & {1} \\
    {1} & {2}   \\
       \end{array} \right),\,
s\left( \begin{array}{cc}
    {0} & {1} \\
    {1} & {0}   \\
       \end{array} \right), \,
(t+1)\left( \begin{array}{cc}
    {0} & {2} \\
    {2} & {0}   \\
       \end{array} \right), \,
(p-2)\left( \begin{array}{c}
    {2}   \\
       \end{array} \right),\,
m\left( \begin{array}{c}
    {0}   \\
       \end{array} \right)).
$$
In particular, the matrix $ A $ has a normal form as in \eqref{eq:normal} where the parameter $ p $ is less than or equal to $ 2 $.

\end{lemma}

\proof
Let $ \mathcal{B} =\{.., g_1,g_2,g_3,...,g_p, ...\}$ be a $ \ZZ $-basis whose Gram matrix is $ A $ such that $ g_1,g_2,...,g_p $ are basis vectors that correspond to $ 1\times 1 $ diagonal blocks $ 2,2,...,2 $.
Then $ (g_i,g_i)\equiv 2 $ for all $ i $ and $ (g_i,g_j)\equiv 0 $, for all $ i \ne j $. 
Let $ \mathcal{B'} $ be the basis obtained from $ \mathcal{B} $ by replacing $ g_1,g_2,g_3 $ by $ g'_1=g_1+g_2, g'_2=g_2+g_3, g'_3=g_1+g_2+g_3 $ respectively. 
Then the change of basis matrix has determinant equal to one, so $ \mathcal{B'} $ is a $ \ZZ $-basis and its Gram matrix modulo $ 4 $ is as in the statement of the lemma 
(i.e. $ (g'_1,g'_1)\equiv (g'_2,g'_2)\equiv 0, (g'_1,g'_2)\equiv 2 $, $ (g'_3,g'_3) \equiv 2$ and  $ (g'_1,g'_3)\equiv (g'_2,g'_3)\equiv 0 \pmod{4} $ and the other products are the same as those for $ \mathcal{B} $). 
This completes the proof of the lemma.

We now show that the parameter $ r $ also has an analogous property:
\begin{lemma}\label{lem:r}
Suppose $ A $ is a matrix such that $ A \pmod{4}$ is the matrix
$$\mathrm{diag}(r
\left( \begin{array}{cc}
    {2} & {1} \\
    {1} & {2}   \\
       \end{array} \right),\,
s\left( \begin{array}{cc}
    {0} & {1} \\
    {1} & {0}   \\
       \end{array} \right), \,
t\left( \begin{array}{cc}
    {0} & {2} \\
    {2} & {0}   \\
       \end{array} \right), \,
p\left( \begin{array}{c}
    {2}   \\
       \end{array} \right),\,
m\left( \begin{array}{c}
    {0}   \\
       \end{array} \right))
$$
for some non-negative integers $p,s,t,m$ and $ r \geq 2 $.
Then, there is a matrix $ P $ with $\det(P)=\mp 1$ such that $P^TAP \pmod{4}$  is 
$$\mathrm{diag}((r-2)
\left( \begin{array}{cc}
    {2} & {1} \\
    {1} & {2}   \\
       \end{array} \right),\,
(s+2)\left( \begin{array}{cc}
    {0} & {1} \\
    {1} & {0}   \\
       \end{array} \right), \,
t\left( \begin{array}{cc}
    {0} & {2} \\
    {2} & {0}   \\
       \end{array} \right), \,
p\left( \begin{array}{c}
    {2}   \\
       \end{array} \right),\,
m\left( \begin{array}{c}
    {0}   \\
       \end{array} \right))
.$$
In particular, the matrix $ A $ has a normal form as in \eqref{eq:normal} where the parameter $ r $ is equal to $ 0 $ or $ 1 $.

\end{lemma}

\proof
Let $ \mathcal{B} =\{e_1,f_1,e_2,f_2...,\}$ be a basis whose Gram matrix is $ A $ such that  $ e_1,f_1,e_2,f_2 $ corresponds to the first two blocks; i.e. for $ i=1,2 $, we have
$ (e_i,f_i)\equiv 1, (e_i,e_i)\equiv (f_i,f_i)\equiv 2\pmod{4}$ and, for $ i \ne j $, $  (e_i,f_j)\equiv 0 \pmod{4}$.
Let 
$ \mathcal{B'} $ be the basis obtained from $ \mathcal{B} $ by replacing $ e_1,f_1,e_2,f_2 $ by $ e'_1=e_1+f_2, f'_1=f_1+e_1+f_2, e'_2=(e_2-f_1)+2(e_1+f_2), f'_2=f_2+2f_1+((e_2-f_1)+2(e_1+f_2)) $ respectively. 
Then the change of basis matrix has determinant equal to one, so $ \mathcal{B'} $ is a $ \ZZ $-basis and, for $ i=1,2 $, we have$ (e'_i,f'_i)\equiv 1, (e'_i,e'_i)\equiv (f'_i,f'_i)\equiv 0 \pmod{4}$ and, for $ i \ne j $, $  (e'_i,f'_j)\equiv 0 \pmod{4}$, so 
$ \mathcal{B'} $ is a $ \ZZ $-basis whose Gram matrix modulo $ 4 $ is as in the statement of the lemma. This completes the proof of the lemma

We now prove a stronger version of Lemma \ref{lem:r} in the case when $ p \geq 1 $:

\begin{lemma}\label{lem:p-r}
Suppose $ A $ is a matrix such that $ A \pmod{4}$ is the following matrix: 
$$\mathrm{diag}(r
\left( \begin{array}{cc}
    {2} & {1} \\
    {1} & {2}   \\
       \end{array} \right),\,
s\left( \begin{array}{cc}
    {0} & {1} \\
    {1} & {0}   \\
       \end{array} \right), \,
t\left( \begin{array}{cc}
    {0} & {2} \\
    {2} & {0}   \\
       \end{array} \right), \,
p\left( \begin{array}{c}
    {2}   \\
       \end{array} \right),\,
m\left( \begin{array}{c}
    {0}   \\
       \end{array} \right))
$$
for some non-negative integers $r,s,t,m$ and $ p \geq 1 $.
Then, for any $ 0\leq r' \leq r+s$, there is a matrix $ P $ with $\det(P)=\mp 1$ such that $P^TAP \pmod{4}$ is 
$$\mathrm{diag}(r'
\left( \begin{array}{cc}
    {2} & {1} \\
    {1} & {2}   \\
       \end{array} \right),\,
(r+s-r')\left( \begin{array}{cc}
    {0} & {1} \\
    {1} & {0}   \\
       \end{array} \right), \,
t\left( \begin{array}{cc}
    {0} & {2} \\
    {2} & {0}   \\
       \end{array} \right), \,
p\left( \begin{array}{c}
    {2}   \\
       \end{array} \right),\,
m\left( \begin{array}{c}
    {0}   \\
       \end{array} \right)).
$$
In particular, the matrix $ A $ has a normal form as in \eqref{eq:normal} where the parameter $ r $ is equal to $ 0 $. 

\end{lemma}

\proof
To prove the lemma, it is enough to show it for $ r'=r-1 $ if $ r\geq 1 $; then, by induction, any matrix as in the hypothesis become congruent to one as in the conclusion with $ r'=0 $.
To proceed, let us assume that $ r\geq 1 $. Let $ \mathcal{B} =\{e_1,f_1,..,g_1,...,g_p,...\}$ be a basis whose Gram matrix is $ A $ such that  $ \{e_1,f_1\} $ corresponds to the first block (for $ r\geq 1 $), and $ g_1,...,g_p $ correspond to the $ 1 \times 1 $ blocks $ 2,2,2,,. $,  so for $ i=1,2 $, we have $ (e_1,f_1)\equiv 1, (e_1,e_1)\equiv (f_1,f_1)\equiv 2\pmod{4}$ and, for $ i \ne j $, $  (e_i,f_j)\equiv 0 \pmod{4}$. Also note that  $ (g_1,g_1)\equiv 2, (g_1,e_1)\equiv (g_1,f_1)\equiv 0 \pmod{4} $.
Let $ \mathcal{B'} $ be the basis obtained from $ \mathcal{B} $ by replacing $ e_1,f_1,g_1 $ by $ e'_1=e_1+g_1, f'_1=f_1+(e_1+g_1), g'_1= g_1+2f_1 $ respectively. 
Then the change of basis matrix has determinant equal to one, so $ \mathcal{B'} $ is a $ \ZZ $-basis and, for $ i=1,2 $, we have $ (e'_1,f'_1)\equiv 1, (e'_1,e'_1)\equiv (f'_1,f'_1)\equiv 0 \pmod{4}, (g'_1,g'_1)\equiv 2, (g'_1,e'_1)\equiv (g'_1,f'_1)\equiv 0 $ and the other products are the same as those for $ \mathcal{B} $. Thus $ \mathcal{B'} $ is a $ \ZZ $-basis whose Gram matrix modulo $ 4 $ is as in the statement of the lemma. 
This completes the proof of the lemma.

\begin{lemma}
\label{lem:q-congruence}
Suppose that $ A $ and $ A' $ are symmetric matrices whose diagonal elements are even numbers. 
Let $ q $ and $ q' $ be the associated quadratic forms (as defined in Sec \ref{sec:intro}). Suppose also that there is a matrix $ P $ with $det(P)=\mp 1$ such that $A'=P^TAP$. 
Then $ q $ and $ q' $ are isomorphic. 


\end{lemma}

\proof 
We denote the associated symmetric forms by $(\, , \, )$ and $(\, , \,)'$ respectively. Then, for any $v $ in $\ZZ^n$, we have $$(Pv,Pv)=(Pv)^TA(Pv)=v^TP^TAPv=v^TA'v=(v,v)'. $$
On the other hand, the quadratic forms $ q $ and $ q' $ are defined as $ q(v)= (v,v)/2 \pmod{2} $ and $ q'(v)= (v,v)'/2 \pmod{2} $ (where, by an abuse of notation, we denote the coset $ v+2\ZZ $ by $ v $). Thus $q(Pv)=q'(v)$, for all $v$ in $\ZZ^n$, i.e. $q$ and $q'$ are isomorphic. This completes the proof of the lemma.

We now find bases for the subspaces defined in \eqref{eq:V000}:
\begin{lemma}
\label{lem:q-V0}

In the setup of Theorem~\ref{th:mod4}, let
$$ \mathcal{B} =\{e_1,f_1,..,,e_r,f_r, ..., e_{r+s}, f_{r+s},g_1,...,g_{2t}, a_1,...,a_p, h_1,...,h_m\}$$ be a
$ \ZZ $-basis whose Gram matrix modulo $ 4 $ is the matrix
$$\mathrm{diag}(r
\left( \begin{array}{cc}
    {2} & {1} \\
    {1} & {2}   \\
       \end{array} \right),\,
s\left( \begin{array}{cc}
    {0} & {1} \\
    {1} & {0}   \\
       \end{array} \right), \,
t\left( \begin{array}{cc}
    {0} & {2} \\
    {2} & {0}   \\
       \end{array} \right), \,
p\left( \begin{array}{c}
    {2}   \\
       \end{array} \right),\,
m\left( \begin{array}{c}
    {0}   \\
       \end{array} \right))
$$
(as in \eqref{eq:normal}). Then the subspaces in \eqref{eq:V000} have the following bases over $ \FF_2 $: 


\vspace{0.1cm}

\noindent
$ \{g_1,...,g_{2t}, a_1,...,a_p, h_1,...,h_m\}$ is a basis for $ \overline{V}_{0} $ and $ \{h_1,...,h_m\}$ is a basis for $ \overline{V}_{000} $ (where, by an abuse of notation, we write $ v $ to denote its coset $ v+2\ZZ $). 

\vspace{0.1cm}

Furthermore:

\vspace{0.1cm}

\noindent
{\rm (i)} if $ p=0 $, then $ \overline{V}_{0} =\overline{V}_{00}$;

\noindent
{\rm (ii)} if $ p=1 $, then $ \{g_1,...,g_{2t}, h_1,...,h_m\}$ is a basis for $ \overline{V}_{00} $; 

\noindent
{\rm (iii)} if $ p=2 $, then $ \{g_1,...,g_{2t},a_1+a_2, h_1,...,h_m\}$ is a basis for $ \overline{V}_{00} $.

\vspace{0.1cm}

In particular, we have the following (where $ q $ is the associated quadratic form as defined in Sec \ref{sec:intro}):

\vspace{0.1cm}

\noindent
{\rm (iv)} $ \dim( \overline{V}_{00} ) \geq 1$ if and only if $ p \geq 2 $ or $ 2t+m \geq 1 $.

\vspace{0.1cm}

\noindent
{\rm (v)} $\dim{\overline{V}_{0}/\overline{V}_{00}}\leq 1$;  $\dim{\overline{V}_{0}/\overline{V}_{00}}= 1$ if and only if $ p\geq 1 $.

\vspace{0.1cm}

\noindent
{\rm (vi)} $q(\overline{V}_0)=0$ if and only if $ p=0 $; $q(\overline{V}_0)=\FF_2$ if and only if $ p\geq 1 $. 

\end{lemma}

\noindent
The proof of the lemma follows from a direct computation. The proof of part (v) uses Lemma \ref{lem:p} and parts (ii,iii).

We can now give linear algebraic interpretations of the parameters in the normal forms of \eqref{eq:normal}:
\begin{lemma}
\label{lem:mod4-unique}
Suppose that $ A $ is a matrix such that $ A \pmod{4} $ is of the form
$$\mathrm{diag}(r
\left( \begin{array}{cc}
    {2} & {1} \\
    {1} & {2}   \\
       \end{array} \right),\,
s\left( \begin{array}{cc}
    {0} & {1} \\
    {1} & {0}   \\
       \end{array} \right), \,
t\left( \begin{array}{cc}
    {0} & {2} \\
    {2} & {0}   \\
       \end{array} \right), \,
p\left( \begin{array}{c}
    {2}   \\
       \end{array} \right),\,
m\left( \begin{array}{c}
    {0}   \\
       \end{array} \right))
$$
(as in \eqref{eq:normal}). Then we have the following:

\noindent
{\rm (i)} $ m= \dim( \overline{V}_{000} ) $,

\noindent
{\rm (ii)} $ 2t+p+m=\dim (\overline{V}_{0}) $,

\noindent
{\rm (iii)} $ p \equiv  \dim(\overline{V}_{0}/\overline{V}_{000}) \pmod{2}$,

\noindent
{\rm (iv)} $ r\equiv Arf(q) \pmod{ 2}$ (if $q(\overline{V}_0)=0$).

\end{lemma}

\proof
The statements (i) and (ii) immediately follow from Lemma \ref{lem:q-V0}. 
For part (iii), we note that $ \dim(\overline{V}_{0}/\overline{V}_{000})= \dim(\overline{V}_{0}) - \dim(\overline{V}_{000})=(2t+p+m)-m=2t+p \equiv  \pmod{2}$, so $ p \equiv  \dim(\overline{V}_{0}/\overline{V}_{000}) \pmod{2}$.

For part (iv), let us recall that the quadratic form $ q $ is defined as $ q(v)\equiv (v,v)/2 \pmod{2}$ (Section \ref{sec:intro}). Then we have the following:

$ q(v)\equiv 1 \pmod{2} $ if and only if $ (v,v) \equiv 2 \pmod{4}  $ and $ q(v)\equiv 0 \pmod{2} $ if and only if $ (v,v) \equiv 0 \pmod{4}  $.

Then, in the setup of Lemma \ref{lem:q-V0}, we have the following: for $ i=1,...,r $, we have $ q(e_i)\equiv q(f_i)\equiv 1 \pmod{2} $ and, for $ i=r+1,...,r+s $, we have $ q(e_i)\equiv q(f_i)\equiv 0 \pmod{2} $. Thus  $q(e_i)q(f_i)\equiv 1 \pmod{2}$ if and only if $ i \in \{1,...,r\} $; taking the sum of these products, we have $ Arf(q) \equiv r\cdot 1 \equiv r \pmod{2}$. This completes the proof of the lemma.

We can now determine normal forms which are congruent modulo $ 4 $ as follows:
\begin{lemma}
\label{lem:normal-conj}
Suppose that $ A $ is a matrix such that $ A \pmod{4} $ is 
$$\mathrm{diag}(r
\left( \begin{array}{cc}
    {2} & {1} \\
    {1} & {2}   \\
       \end{array} \right),\,
s\left( \begin{array}{cc}
    {0} & {1} \\
    {1} & {0}   \\
       \end{array} \right), \,
t\left( \begin{array}{cc}
    {0} & {2} \\
    {2} & {0}   \\
       \end{array} \right), \,
p\left( \begin{array}{c}
    {2}   \\
       \end{array} \right),\,
m\left( \begin{array}{c}
    {0}   \\
       \end{array} \right))
$$

and $ A' $ is a matrix such that $ A' \pmod{4} $ is 
$$\mathrm{diag}(r'
\left( \begin{array}{cc}
    {2} & {1} \\
    {1} & {2}   \\
       \end{array} \right),\,
s'\left( \begin{array}{cc}
    {0} & {1} \\
    {1} & {0}   \\
       \end{array} \right), \,
t'\left( \begin{array}{cc}
    {0} & {2} \\
    {2} & {0}   \\
       \end{array} \right), \,
p'\left( \begin{array}{c}
    {2}   \\
       \end{array} \right),\,
m'\left( \begin{array}{c}
    {0}   \\
       \end{array} \right))
$$
(as in \eqref{eq:normal}). 
Let $ q $ and $ q' $ be the quadratic forms associated with $ A $ and $ A' $ respectively.
Then there is a matrix $ P $ with $ \det(P)=\pm 1 $ such that $ P^TAP\equiv A' \pmod{4} $ if and only if $ m=m' $ and the quadratic forms $ q $ and $ q' $ are isomorphic.

\end{lemma}

\proof
Let us suppose that there is a matrix $ P $ with $ \det(P)=\pm 1 $ such that $ P^TAP\equiv A' \pmod{4}$. Then 
the subspaces  $ \overline{V}_{000} $ and the quadratic form $ q $ are defined in terms of the bilinear form $ (\, ,\, ) $ associated to $ A $, so the dimension of $ \dim{\overline{V}_{000} }$ and the isomorphism class of $ q $ are invariant under the matrix congruence relation (see also Lemma \ref{lem:q-congruence}). This shows that $ m=m' $ and the forms $ q,q' $ are isomorphic.

Conversely, let us suppose that $ m=m' $ and $ q,q' $ are isomorphic. We denote the kernels of the associated bilinear forms by $ \overline{V}_{0} $ and $ \overline{V}'_{0} $. Since $ q $ and $ q' $ are isomorphic, we have $ \dim(\overline{V}_{0})=\dim(\overline{V}'_{0}) $ and $ q(\overline{V}_{0})= q'(\overline{V}'_{0})$ . Also, according to the classification of the quadratic forms, we have the following:

if $ q(\overline{V}_{0})= q'(\overline{V}'_{0})=0$, then $ Arf(q)=Arf(q') $. Also, by Lemma \ref{lem:mod4-unique}, $ r\equiv Arf(q) \pmod{ 2} $ and $ r'\equiv Arf(q') \pmod{ 2} $, so $ r\equiv r' \pmod{ 2} $. Then, by Lemma \ref{lem:r}, there is a matrix $ P $ with $ \det(P)=\pm 1 $ such that $ P^TAP\equiv A' \pmod{4}$.

if $ q(\overline{V}_{0})= q'(\overline{V}_{0})=\FF_2$, then $ p,p' \geq 1$ (Lemma \ref{lem:q-V0}). Also, by Lemma \ref{lem:mod4-unique}, $ p \equiv  \dim(\overline{V}_{0}/\overline{V}_{000})= \dim(\overline{V}_{0})- \dim(\overline{V}_{000}) \pmod{ 2}$. Also $ m'=m= \dim(\overline{V}_{000})$ and, since $ q $ and $ q' $ are isomorphic, $ \dim(\overline{V}_{0})=\dim(\overline{V}'_{0}) $, so $ p\equiv p' \pmod{ 2} $. Then, by Lemma \ref{lem:p}, there is a matrix $ P $ with $ \det(P)=\pm 1 $ such that $ P^TAP\equiv A' \pmod{4}$.
This completes the proof of the lemma.

Let us now prove Theorem~\ref{th:mod4-unique}. To prove the theorem, we first note that the subspace $ \overline{V}_{000} $ and the quadratic form $ q $ are defined in terms of the bilinear form $ (\, ,\, ) $ 
associated to $ A $, so the dimension $ \dim{\overline{V}_{000} }$ and the isomorphism class of $ q $ are invariant under the matrix congruence relation (see also Lemma \ref{lem:q-congruence}).
Thus, to prove Theorem~\ref{th:mod4-unique}, it is enough to assume that $ A \pmod{4} $ is as in \eqref{eq:normal}. Then, by Lemma \ref{lem:normal-conj}, the congruence class of $ A $ is determined by $ \dim{\overline{V}_{000} }$ and the isomorphism class of $ q $. The uniqueness part of the theorem also follows from Lemma \ref{lem:normal-conj}. More explicitly, if $q(\overline{V}_0)=0$ then $ p=0 $ (by Lemma \ref{lem:q-V0}(vi)) and, by Lemma \ref{lem:r}, we can take $ r \in \{0,1\} $. Similarly, if $q(\overline{V}_0)=\FF_2$, then $ p\geq 1 $, so, by Lemma \ref{lem:p}, we can take $ p \in \{1,2\} $ and, by Lemma \ref{lem:p-r}, $r=0 $. This completes the proof of Theorem~\ref{th:mod4-unique}.






\subsection{Proof of Corollary~\ref{cor:det}} By Theorem \ref{lem:mod4}, we can assume that $ A \pmod{4}$ is of the form 
$$\mathrm{diag}(r
\left( \begin{array}{cc}
    {2} & {1} \\
    {1} & {2}   \\
       \end{array} \right),\,
s\left( \begin{array}{cc}
    {0} & {1} \\
    {1} & {0}   \\
       \end{array} \right), \,
t\left( \begin{array}{cc}
    {0} & {2} \\
    {2} & {0}   \\
       \end{array} \right), \,
p\left( \begin{array}{c}
    {2}   \\
       \end{array} \right),\,
m\left( \begin{array}{c}
    {0}   \\
       \end{array} \right)),
$$ where $ r,s,p,t,m $ are non-negative integers and $ n= 2(r+s)+2t+m$; also $ \dim(\overline{V}_{0})=2t+p+m$. Then $ \det{A}=(-1)^{r+s}\cdot 2^p\cdot  0^{2t+m} \pmod{4} $, as a product of the determinants of the blocks, so we have the following:

\begin{itemize} 
\item
$  \det{A}\equiv (-1)^{r+s} \pmod{4} $ if and only if $ {2t+m}= 0 , p=0$. Then $ r+s $ is even (resp. odd) if and only if $ n= 2(r+s)+p+2t+m=2(r+s)\equiv 0 \pmod{4} $ (resp. $ n \equiv 2 \pmod{4}  $ ). Thus $ \det{A}$ is as stated in the corollary.
\item
$  \det{A}\equiv 2 \pmod{4} $ if and only if $ 2t+m = 0$ and $ p=1 $. This happens, by Lemma \ref{lem:q-V0}, if and only if $ \dim(\overline{V}_{00}) = 0  $ and $ \dim(\overline{V}_{0}) = 1  $.
\item

$  \det{A}\equiv 0 \pmod{4} $ if and only if $ 2t+m \geq 1$ or $ p \geq 2 $. This happens, by Lemma \ref{lem:q-V0}, if and only if $ \dim(\overline{V}_{00}) \geq 1$.

\end{itemize} 







It follows, in particular, that if the rank of the matrix $ A \pmod{2} $ is less than or equal to $ n-2 $ (i.e. $ \dim(\overline{V}_{0}) \geq 2  $), then $ \dim(\overline{V}_{00})\geq \dim(\overline{V}_{0})-1 \geq 1$ (by Lemma \ref{lem:q-V0}), so $  \det{A}\equiv 0 \pmod{4} $. 
This completes the proof of the corollary.

\subsection{Proof of Theorem~\ref{th:mut}} \label{sec:mut-proof}
Let $ q $ and $ q' $ be the associated quadratic forms for $ A $ and $ A' $ respectively. 
Then, by Corollary \ref{cor:skew-quad}, the forms $ q $ and $ q' $ are isomorphic. Thus, by Corollary \ref{cor:det} and Remark \ref{rem:det}, we have $  \det{A}\equiv \det{A'}  \pmod{4} $. 
Also, by Theorem~\ref{th:mod4-unique}, the normal forms of $ A $ and $ A' $ are determined by the associated quadratic forms $ q $ and $ q' $ with the parameters $ m $ and $ m' $ respectively. This implies that, since $ q $ and $ q' $ are isomorphic, the matrix-congruence of $ A $ and $ A' $ modulo $ 4 $ is determined by $ m $ and $ m' $ being equal. This completes the proof of the theorem.


To give an example, the following matrices are quasi-Cartan companions of the same matrix but the corresponding parameters are $ m=0 $ and $ m'=2 $ respectively:

$$A=\left(\begin{array}{rrrr}
2 & 1 & 0 & 1  \\
1 & 2 & 1 & 0  \\
0 & 1 & 2 & -1 \\
1 & 0 & -1 & 2
\end{array}\right),\qquad
A'=\left(\begin{array}{rrrr}
2 & 1 & 0 & 1  \\
1 & 2 & 1 & 0  \\
0 & 1 & 2 & 1 \\
1 & 0 & 1 & 2
\end{array}\right).$$
In fact, for the matrix 
$$P=\left(\begin{array}{rrrr}
1 & 0 &-1&2   \\
0 & 1 & 2  &-1 \\
0 & 0 & 1 &0\\
0 & 0 & 0 &1
\end{array}\right),$$ we have the following normal forms  $ P^{T}AP $ and $ P^{T}A'P $ $ \pmod 4 $ of $ A $ and $ A' $ respectively:

$$\left(\begin{array}{rrrr}
2 & 1 & 0 &0  \\
1 & 2 & 0 & 0  \\
0 & 0 & 0 & 2 \\
0 & 0 & 2 & 0
\end{array}\right),\qquad
\left(\begin{array}{rrrr}
2 & 1 & 0 &0  \\
1 & 2 & 0 & 0  \\
0 & 0 & 0 & 0 \\
0 & 0 & 0 & 0
\end{array}\right).$$

\subsection{Example for  Corollary~\ref{cor:delta}} \label{sec:delta-ex}
We consider the example given in  \cite[Proposition 1.7]{C} for $ n=3 $. More precisely, we consider the following skew-symmetric matrices: 

$$B=\left(\begin{array}{rrr}
0 & 1 & 0   \\
-1 & 0 & 1   \\
0 & -1 & 0
\end{array}\right),\qquad
B'=\left(\begin{array}{rrr}
0 & 2 & 0   \\
-2 & 0 & 1   \\
0 & -1 & 0
\end{array}\right).$$
As noted in \cite[Section 3.2]{C}, the matrices $ B $ and $ B' $ are congruent over $\ZZ $: for the matrix 
$$P=\left(\begin{array}{rrr}
1 & 0 &0   \\
0 & 1 & 0   \\
-1 & 0 & 1
\end{array}\right),$$ we have $ P^{T}BP=B' $. As in Corollary~\ref{cor:delta}, let $ \overline{B} $ (resp. $ \overline{B'} $) denote the matrix $ B \pmod{2} $ (resp. $ B' \pmod{2} $). Then both $  \overline{B} $ and $  \overline{B'} $ have rank equal to two: their kernels (i.e. $ \overline{V}_{0} $ and $\overline{V}_{0}')$ are spanned by the vectors $ u=(1,0,1) $ and $ u'=(1,0,0) $ respectively. Also, for the corresponding quadratic forms $ Q=Q_{B} $ and $ Q'=Q_{B'} $, we have the following: $ Q(u)=0 $ and $ Q'(u')=1 $. Thus,  in the setup of Corollary~\ref{cor:delta}, we have $ Q_B(\overline{V}_{0})=\{0\}$ and $ Q_{B'}(\overline{V}'_{0}) = \FF_2$ (in particular, $ Q_B $ and $ Q_{B'} $ are not isomorphic). Then, by the corollary, $ \delta(B) \ne \delta(B') $. In fact, more generally (by Corollary \ref{cor:det}), for any quasi-Cartan companions $ A $ and $ A' $ of $ B $ and $ B' $ respectively, we have $ \det(A)\equiv 0 \pmod 4 $ and $ \det(A') \equiv 2 \pmod 4 $. (In particular, $ \delta(B)=0 $ and $ \delta(B')=2 $). This can be checked, for example, on the following particular quasi-Cartan companions:

$$A=\left(\begin{array}{rrr}
2 & 1 & 0   \\
1 & 2 & 1   \\
0 & 1 & 2
\end{array}\right),\qquad
A'=\left(\begin{array}{rrr}
2 & 2 & 0   \\
2 & 2 & 1   \\
0 & 1 & 2
\end{array}\right).$$
A simple calculation shows indeed that $ \det(A)\equiv 0 \pmod 4 $ and $ \det(A') \equiv 2 \pmod 4 $. Let us also note that, since the quadratic forms $ Q_B $ and $ Q_{B'} $ are not isomorphic, the matrices $ B $ and $ B' $ are not mutation equivalent by Proposition \ref{prop:skew-quad}.

\end{document}